\documentclass[12pt, reqno, a4paper]{amsart}

\usepackage{amsthm,amssymb,amsmath,amsfonts}
\usepackage{verbatim}
\usepackage{bbm}
\usepackage{hyperref}
\usepackage[all]{xy}
\usepackage{latexsym}
\usepackage{color}
\usepackage{tikz, tikz-qtree}
\usepackage{graphicx}
\usepackage{caption, subcaption}
\usepackage{verbatim}
\usepackage{rotating} 
\usepackage[colorinlistoftodos]{todonotes} 

\usepackage{amsaddr}

\usepackage{natbib}
\usepackage{setspace}
\newtheorem{thm}{Theorem}
\newtheorem{prop}[thm]{Proposition}

\theoremstyle{definition}

\def\R{\mathbb{R}}

\begin{document}

\title{From Pareto to Weibull - a constructive review of distributions on $\R^+$}
\author{Corinne Sinner$^*$,  Yves Dominicy$^\dagger$,  Julien~Trufin$^*$, Wout Waterschoot$^\ddagger$, Patrick Weber$^*$ and Christophe Ley$^{\S\ddagger}$}
\address{$^*$D\'epartement de Math\'ematiques, Universit\'e libre de Bruxelles CP210, Boulevard du Triomphe, Bruxelles 1050, Belgique.}
\address{$^\dagger$ Bank employee, Luxembourg.}
\address{$^\ddagger$Department of Applied Mathematics, Computer Science and Statistics, Ghent University, Krijgslaan 281 - S9, 9000 Gent, Belgium}
\address{$^\S$Department of Mathematics, University of Luxembourg, 6, rue de la Fonte, 4365 Esch-sur-Alzette, Luxembourg}

%

\maketitle

\begin{abstract}
Power laws and power laws with exponential cut-off are two distinct families of  distributions on the positive real half-line. In the present paper, we propose a unified treatment of both families by building a family of distributions that interpolates between them, which we call Interpolating Family (IF) of distributions. Our original construction, which relies on techniques from statistical physics,  provides a connection for hitherto unrelated distributions like the Pareto and Weibull distributions, and sheds new light on them. The IF also contains several distributions that are neither of power law nor of power law with exponential cut-off type. We calculate quantile-based properties, moments and modes  for the IF. This allows us to review known properties of famous  distributions on $\R^+$ and to provide in a single sweep these characteristics for various less known (and new) special cases of our Interpolating Family. 
\end{abstract}

\hfill

\noindent%


\section{Introduction}
\label{sec:intro}

Initiated in the 19th century by famous mathematicians as Adolphe Quetelet, Sir Francis Galton or Vilfredo Pareto, the journey for the probability distribution that best describes observations has never ceased since then. Still nowadays it remains among the most popular topics in statistics as shown by the large amount of scientific papers published on the subject (see for instance the review papers by~\citet{Jones2015} and~\citet{Ley2015}). 

In this paper, we are concerned with probability distributions that analyse data on $\R^+$, which we refer to as size-type data. Size distributions are probability laws designed to model data that only take non-negative values or values above a certain threshold. Typical examples of such data are claim sizes in actuarial science, wind speeds in meteorology or lifetime data. Nonetheless, the spectrum of application areas is much broader and non-negative observations appear naturally in survival analysis~\citep{lawless2003, leewang2003}, environmental science~\citep{marchenko2010, White2008}, network traffic modeling~\citep{mitzenmacher2004}, web hits (\citet{Crovella1997}, \citet{huberman1999}  and \citet{Adamic2000}), reliability theory~\citep{rausandhoyland2004}, astrophysics \citep{Aschwanden2013}, economics \citep{eeckhout2004, luttmer2007, Farmer2008, piketty2014, toda2015, gabaix2016}, income of top earners in areas of arts, sports and business~\citep{rosen1981}, seismology~\citep{ley2020modelling, Burroughs2001},  hydrology~\citep{clarke2002, aban2006},  biological systems \citep{Munoz2018}, neuroscience \citep{Chialvo2010}, counting word frequencies (\citet{Estoup1916} and \citet{Zipf1949}), counting  citations of scientific papers \citep{Price1965},  the distribution of the number of calls received on a single day (\citet{Aiello2000} and \citet{Huberman2004} ), the study of diameter of moon craters (\citet{Neukum1994}), the intensity of solar flares \citep{Lu1991}, the intensity of wars (\citet{Small1982} and  \citet{Roberts1998}), counting the frequencies of family names in a country \citep{Miyazima2000}, to cite but these. In geoscience, interest in size distributions appeared in the 80s of the past century. It was Mandelbrot’s work on fractals ~\citep{Mandelbrot1983} which drew attention to the distribution of sizes of diverse geological objects and structures, like lakes, faults, fault gouge, oil reservoirs, sedimentary layers, and even asteroids ~\citep{Turcotte1997}.
Given the range of distinct domains of application, it is not surprising that there exists a plethora of different size distributions and that it is still a very active research area (see for instance \citet{Kleiber2003}, \citet{sornette2003}, \citet{mitzenmacher2004} and \citet{dominicy2017distributions}).

\hfill

Parametric distributions of size phenomena is a subject that has been explored in statistical literature for over 100 years since the publication of the Italian economist and engineer Vilfredo Pareto's famous book ‘Cours déconomie politique’ in 1897 (\citet{Pareto1897}). Vilfredo observed in 1895 that in many populations the number of individuals with income exceeding a given threshold can be approximated by a probability law, nowadays known as a Pareto distribution \citep{Arnold2015}. Pareto used his law to model the distribution of income and the allocation of wealth among individuals. In 1912, the Norwegian actuary Birger Meidell used the Pareto law to model maximum risk in life insurance \citep{Meidell1912}. He based his study on the hypothesis that the insured sum is proportional to the income of the policy holder. Some years later, \citet{Auerbach1916} showed that the population of cities and the frequency of words in texts, respectively, follow essentially the same statistical pattern \citep{Newman2005}. Over the years, the Pareto law has further been applied to city size, file size distribution of internet traffic, the size of meteorites, or the size of sand particles~\citep{reed2004}. \citet{Mandelbrot1964} derived a Pareto distribution of the amount of fire damage from the assumption that the probability of the fire increasing its intensity at any instant of time is constant. For a discussion in detail about the Pareto distribution we refer the reader to the books by \citet{Kleiber2003} and \citet{Arnold2015} .

This very popular size distribution, also called the Pareto type I distribution, has a probability density function
$$x \mapsto  \frac{\alpha x_{\mathrm{0}}^{\alpha}}{x^{\alpha+1}}, \quad x \in [x_0,\infty),$$
where $x_0\geqslant0$ is a location parameter and $\alpha>0$ is a shape parameter known as the tail or Pareto index. Note that a decreasing value of $\alpha$ implies a heavier tail.
The Pareto distribution is a member of the power laws, which are typically of the form $x\mapsto k x^{-\alpha}$, with normalizing constant $k$ and power $\alpha >0$.

\hfill

A popular alternative to power laws are power laws with exponential cut-off, for which the Weibull distribution is a famous representative of, whose densities take the form $x \mapsto k x^{-\alpha}\mathrm{e}^{-\beta x}$ with normalizing constant $k$, power exponent $\alpha >0$ and rate parameter $\beta>0$.
A power law with exponential cut-off behaves like a power law for small values of $x$, while its tail behavior is governed by a decreasing exponential.
The Weibull distribution has density
$$x \mapsto \frac{\alpha}{\sigma}\left(\frac{x}{\sigma}\right)^{\alpha-1}e^{-\left(\frac{x}{\sigma}\right)^{\alpha}}, \quad x \in [0,\infty),$$
where $\alpha > 0$ is a shape parameter regulating tail-weight and $\sigma > 0$ is a scale parameter.

Pareto's contribution simulated further research in the specification of new models to fit the whole range of income. In 1898, the French statistician Lucien March proposed to use the gamma distribution to fit the distribution of wages in France, Germany, and the USA.  Note that the Weibull law is a generalized gamma distribution which is a generalization of the gamma distribution.

Historically, the Weibull distribution it gained its prominence and name by the Swedish engineer and scientist Waloddi Weibull, who discussed it in his paper \citet{Weibull1939}. He described the distribution in detail in his work \citet{Weibull1951}, although it was actually first identified by \citet{Frechet1927} and already applied by \citet{Rosin1933} to describe a particle size distribution.
Nowadays, the Weibull distribution is widely used in various domains such as life data analysis~\citep{nelson2005}, wind speed modeling~\citep{manwell2009} and hydrology~\citep{clarke2002}.

\hfill

Power laws and power laws with exponential cut-off are mostly studied apart from each other, due to their disparity. In the present paper, we shall build a bridge between these two classes of size distributions by proposing an over-arching family of distributions that interpolates between both classes, hence the name \emph{Interpolating family of size distributions} (for simplicity, we shall from now on also refer to it as IF distribution). The corresponding density, which we shall discuss in detail later in the paper, is of the form
\begin{equation}\label{newdist}
        \text{IF}(x;p,b,c,q,x_0)=\frac{|b|q}{c}\left(\frac{x-x_0}{c}\right)^{b-1} G_p(x)^{-q-1}\left(1- \frac{1}{p+1}G_p(x)^{-q}\right)^p,
\end{equation}
for $x\in [x_0,\infty)$, and with $p\in[0,\infty], b\neq0, c,q>0, x_0\geqslant 0$ and
\begin{equation*}
        G_p(x)=(p+1)^{-\frac{1}{q}}+ \left(\frac{x-x_0}{c}\right)^b.
\end{equation*}

The roles of the various parameters will be described in Section~\ref{sec:para}. The alert reader will no doubt have recognized the densities of various size distributions included in~\eqref{newdist}, such as the Pareto and Weibull. The essence of the IF  rests on its construction: we have used a technique from statistical mechanics (see Section~\ref{sec:IF}) that allows us to interpolate between the Pareto and Weibull distributions, even more generally, between power laws and power laws with exponential cut-off. Thus, we are precisely  finding a path from one end of the spectrum to the other, and this moreover in a constructive way. Besides providing a link between these \emph{a priori} distinct families of size distributions, our proposal also permits to treat their properties such as moments, quantiles, modes in a unique way. Thus, by studying these properties of the IF, we are reviewing existing characteristics for certain famous distributions and at the same time we are uncovering these properties for less studied size distributions.

\begin{figure}
		\includegraphics[scale=0.6]{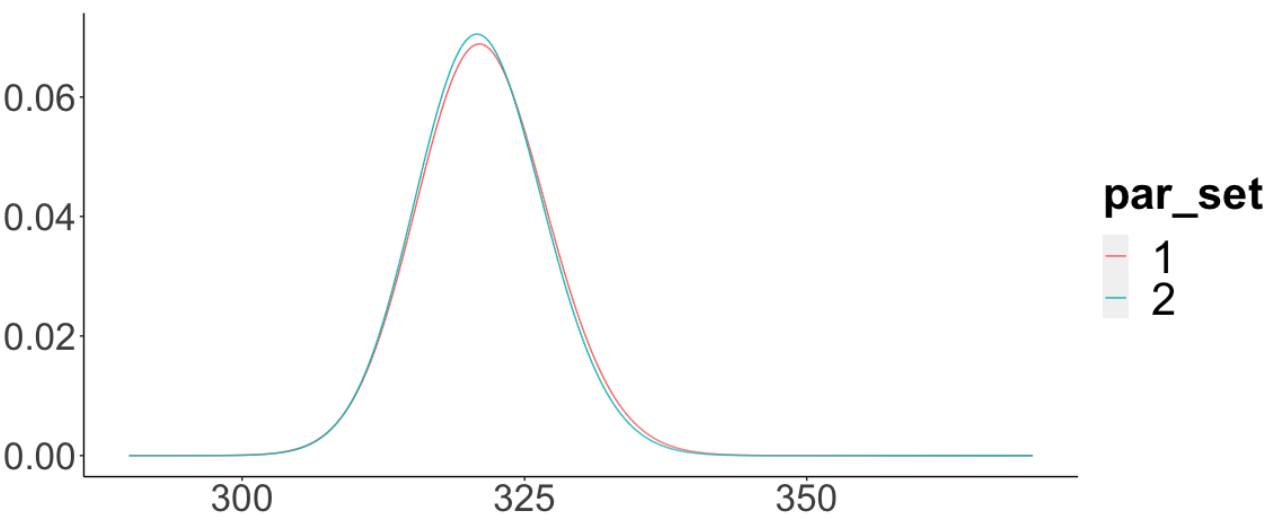}
		\caption{IF density for two parameter sets: set 1 (red line) corresponds to $p = 10, b = 20, c = 300, q = 20, x_0 = 50$ and set 2 (blue line) to $p = 11.33, b = 20.22, c = 322.40, q = 34.08, x_0 = 35.98$.}
		\label{closeness}
	\end{figure}

We wish to stress that the goal of our paper is very different from most flexible modelling papers. With the IF distribution we wish to provide new insights into a common constructive root of apparently disjunct families of distributions, and review their properties from a novel standpoint. Our aim is thus NOT to build a new model that one should use to fit various data sets. Indeed, one has to be cautious with the IF distribution and its 5 parameters. For instance, maximum likelihood estimation of all parameters
involved in (\ref{newdist}) is far from trivial, especially the parameters $x_0$ and $p$ may require a separate treatment in the full 5-parameter model. Numerical optimization methods can often land in different local optima, which is further accentuated by the fact that distinct parameter combinations lead to very similar density shapes; see Figure~\ref{closeness}. Our own investigations seem to imply that the 5-dimensional surface which corresponds to the log-likelihood function behaves very chaotically, which makes the maximization procedure for algorithms which are designed for that purpose extremely difficult. Therefore we do not dwell upon inferential and computational issues, which require a complete treatment on their
own (and for which we have not yet found a satisfying solution, having tried out likelihood-based, moment-based and quantile-based approaches).  This goes beyond the scope of the present
paper, and we do refer the interested reader to \cite{Clauset98} for
an insightful discussion and suggestions for caution (``Commonly used
methods for analyzing power-law data, such as least-squares fitting, can
produce substantially inaccurate estimates of parameters for power-law
distributions...'') regarding parameter estimation for power laws and
power laws with exponential cut-off.
\hfill

The remainder of the paper is organized as follows. In Section~\ref{sec:IF} we describe the Interpolating Family of size distributions, explain how it interpolates between power laws and power laws with exponential cut-off and elucidate the role of each of the five parameters. In Section~\ref{sec:special} we summarize some of the special cases of the IF. We then provide quantile-based properties in Section~\ref{sec:quant}, moment-based results in Section~\ref{sec:moments}, and mode-related results in Section~\ref{sec:mode}. 
Final comments are stated in Section~\ref{sec:final}, while technical derivations are presented in the Appendix.


\section{The Interpolating Family: construction and parameter interpretation} 
\label{sec:IF}

In this section we present in detail our original construction leading to the Interpolating Family of size distributions. Section~\ref{sec:para} expounds on the role of each of the five parameters. 


\subsection{Construction of the family} 
\label{sec:meth}

As announced in the Introduction, our goal is to build a size distribution which incorporates both power laws and power laws with exponential cut-off. 
To show that~\eqref{newdist} indeed satisfies this requirement, we start by writing up power law distributions and  power law distributions with exponential cut-off in a unified language.

\subsubsection{Power laws}

The probability density function (pdf) of a typical power law distribution corresponds to
$$x\mapsto q(1+x)^{-q-1}, \quad x \in [0,\infty),$$
where the tail behavior is governed by the shape parameter $q>0$. To get a more flexible distribution, one may add various parameters, such as a scale parameter $c>0$, a location parameter $x_0 \geqslant 0$ and/or a shape parameter $b>0$, leading to
\begin{equation*}
	x\mapsto \frac{bq}{c} \left(\frac{x-x_0}{c}\right)^{b-1} \left(1+\left(\frac{x-x_0}{c}\right)^b\right)^{-q-1}, \quad x \in [x_0,\infty).
\end{equation*}
Alternatively, in terms of the function $G_0(x)=1+\left(\frac{x-x_0}{c}\right)^b$, the pdf can be written under the form
\begin{equation*}
	f_0(x)=q~g_0(x) G_0(x)^{-q-1},  \quad x \in [x_0,\infty), 
\end{equation*}
where $g_0(x)=\frac{\mathrm{d}}{\mathrm{d}x}G_0(x)=\frac{b}{c}\left(\frac{x-x_0}{c}\right)^{b-1}$. We point out that the function $G_0(x)$ has been chosen such that the following boundary conditions are satisfied: $G_0(x_0)=1$ and ${\lim \limits_{x \to \infty} G_0(x)=\infty}$.

\subsubsection{Power laws with exponential cut-off}

The pdf of a typical power law distribution with exponential cut-off reads
\begin{equation*}
	x\mapsto qx^{-q-1}e^{-x^{-q}}, \quad x \in [0,\infty).
\end{equation*}
The shape parameter $q>0$ still controls the tail behavior and, just as for power laws, we may increase the flexibility of the model by adding scale, location and shape parameters to get
\begin{equation*}
	x\mapsto\frac{bq}{c} \left(\frac{x-x_0}{c}\right)^{-bq-1} e^{-\left(\frac{x-x_0}{c}\right)^{-bq}}, \quad x \in [x_0,\infty).
\end{equation*}
Alternatively, we may write the pdf in terms of the function $G_{\infty}(x)=\left(\frac{x-x_0}{c}\right)^b$ as
\begin{equation*}
	f_{\infty}(x)=q~g_{\infty}(x) G_{\infty}(x)^{-q-1} e^{-G_{\infty}(x)^{-q}}, \quad x \in [x_0,\infty), 
\end{equation*}
where $g_{\infty}(x)=\frac{\mathrm{d}}{\mathrm{d}x}G_{\infty}(x)=\frac{b}{c}\left(\frac{x-x_0}{c}\right)^{b-1}$. Note that the function $G_{\infty}(x)$ has been chosen such that ${G_{\infty}(x_0)=0}$ and $\lim \limits_{x \to \infty} G_{\infty}(x)=\infty$.

\subsubsection{Interpolating Family} 

If we want a highly flexible distribution including both power laws and power laws with exponential cut-off, we need a way to build densities interpolating between  $f_0(x)$ and $f_{\infty}(x)$. To  this end, we introduce a mild variant of the one-parameter deformation of the exponential function popularized in the seminal paper \cite{tsallis1988} in the context of non-extensive statistical mechanics. A more detailed account can be found in  the review paper~\cite{tsallis2002}.

\hfill

For any $p \in [0,\infty]$, we define the $p$-exponential\footnote{The classical $q$-exponential defined by \cite{tsallis1988}  has the form $\tilde{e}_q(x)= \left(1+(1-q)x\right)^{\frac{1}{1-q}}$. 
We slightly modified the deformation path in order to simplify the calculations.} by
\begin{equation*}
	e_p(x)=\left(1-\frac{1}{p+1} x\right)^p, \quad x \in [0,p+1].
\end{equation*}
The extreme cases $p=0$ and  $p \to  \infty$ respectively correspond to  $e_0(x)=1$ over~$[0,1]$ and $e_\infty(x)=e^{-x}$ over $[0,\infty)$. With this in mind, it is natural to consider densities of the type
\begin{equation}\label{firsttry}
	f_p(x)= q~g_p(x) G_p(x)^{-q-1} e_p\left(G_p(x)^{-q}\right), \quad x \in [x_0,\infty), 
\end{equation}
with $g_p(x)=\frac{\mathrm{d}}{\mathrm{d}x}G_p(x)$, where we have not defined the function $G_p(x)$ yet. Just as $e_p(x)$ interpolates between 1 and $e^{-x}$, the mapping $G_p$ should also vary between $G_0$ and $G_\infty$. 
Hence, with the parameters $c>0$, $b>0$ and $x_0 \geqslant 0$ bearing the same interpretation as before, the map $G_p$ could be chosen as $G_p(x)=k+\left(\frac{x-x_0}{c}\right)^b$ for some constant $k$.
A quick calculation shows that $k=(p+1)^{-\frac{1}{q}}$ is the right choice for $f_p(x)$ to integrate to one over its domain. Consequently
\begin{equation*}
	G_p(x)=(p+1)^{-\frac{1}{q}}+ \left(\frac{x-x_0}{c}\right)^b, \quad x \in [x_0, \infty).
\end{equation*}
Since $G_p(x)^{-q}$ with $q>0$  maps $[x_0, \infty)$ onto $[0, p+1]$, the function $e_p(G_p(x)^{-q})$ is well-defined. The pointwise convergence of the resulting density $f_p$ to $f_0$ as $p$ tends to zero (respectively $f_p$ to $f_\infty$ as $p\to \infty$) can be shown by straightforward limit calculations which we omit here. 

\hfill

The density~\eqref{firsttry} now almost corresponds to the density announced in the Introduction. Relaxing the condition $b>0$ into $b\in\R_0$, we finally end up with 
\begin{equation}\label{IF}
\text{IF}(x;p,b,c,q,x_0)=\text{sign}(b)q~g_p(x) G_p(x)^{-q-1} e_p\left(G_p(x)^{-q}\right), \quad x\in [x_0,\infty).
\end{equation}
The relaxation on $b$ only entails a minor change in the normalizing constant, which remains extremely simple. We call IF the \emph{Interpolating Family} of size distributions as it interpolates between power laws and power laws with exponential cut-off. The density depends on five parameters $p,b,c,q$ and $x_0$, which we will discuss in more detail in the next section. 


\subsection{Interpretation of the parameters}\label{sec:para} 

\begin{figure}[h!]
	\centering
	\begin{subfigure}{0.48\textwidth}
		\includegraphics[scale=0.5]{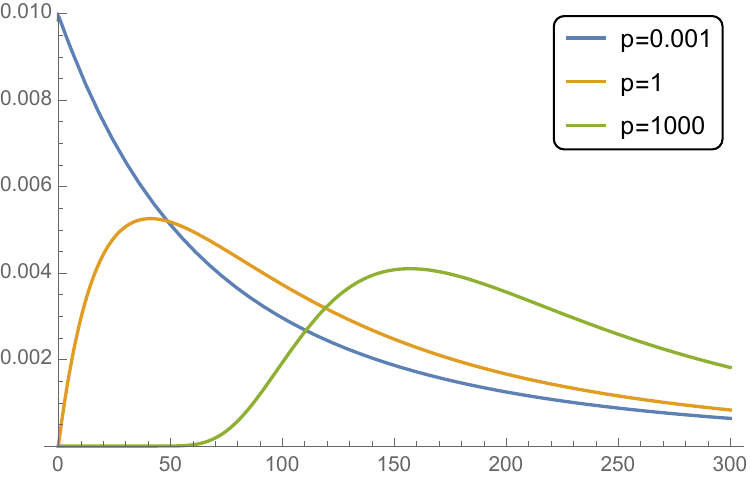}
	\end{subfigure}
	\hfill
	\begin{subfigure}{0.48\textwidth}
		\includegraphics[scale=0.5]{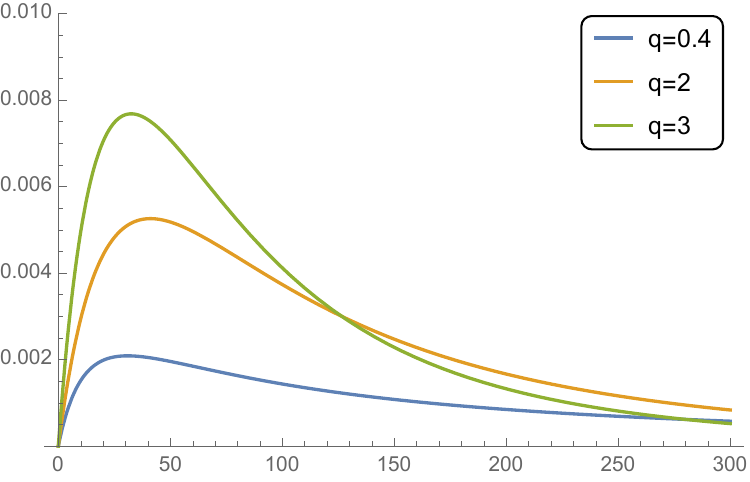}
	\end{subfigure}
	\centering
	\begin{subfigure}{0.48\textwidth}
		\includegraphics[scale=0.5]{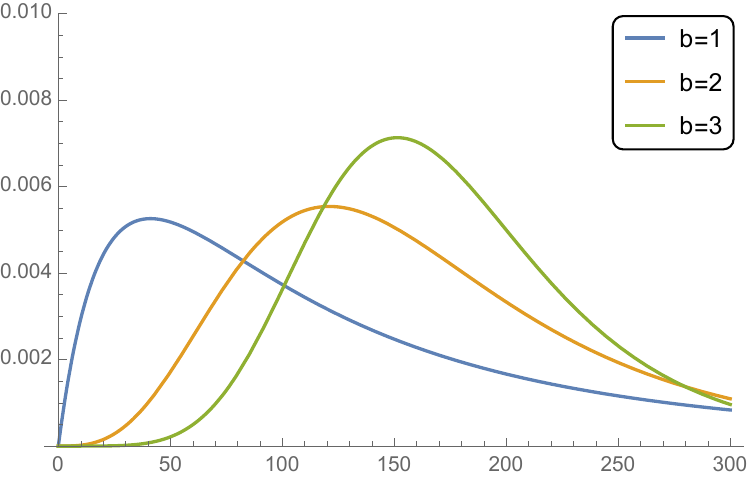}
	\end{subfigure}
	\hfill
	\begin{subfigure}{0.48\textwidth}
		\includegraphics[scale=0.5]{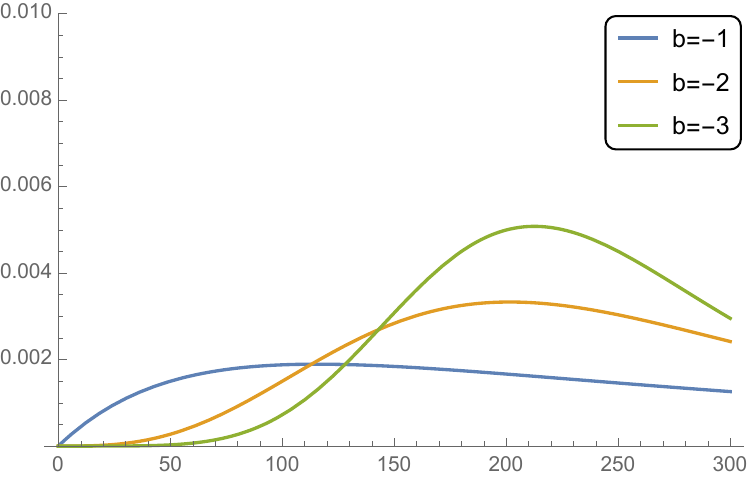}
	\end{subfigure}
	\centering
	\begin{subfigure}{0.48\textwidth}
		\includegraphics[scale=0.5]{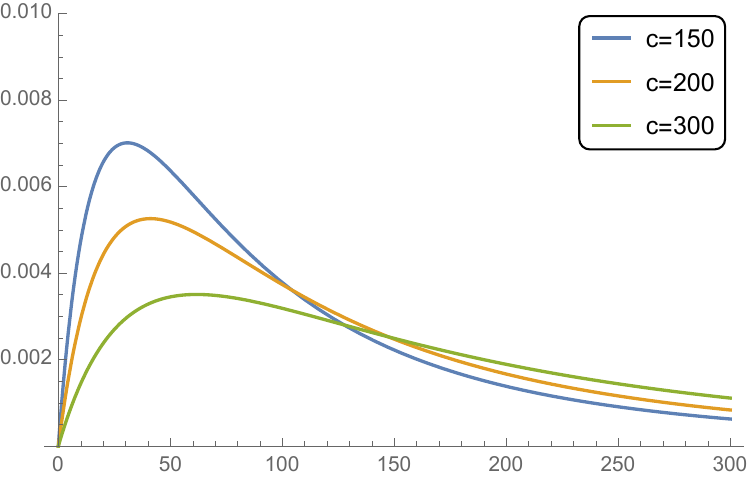}
	\end{subfigure}
	\hfill
	\begin{subfigure}{0.48\textwidth}
		\includegraphics[scale=0.5]{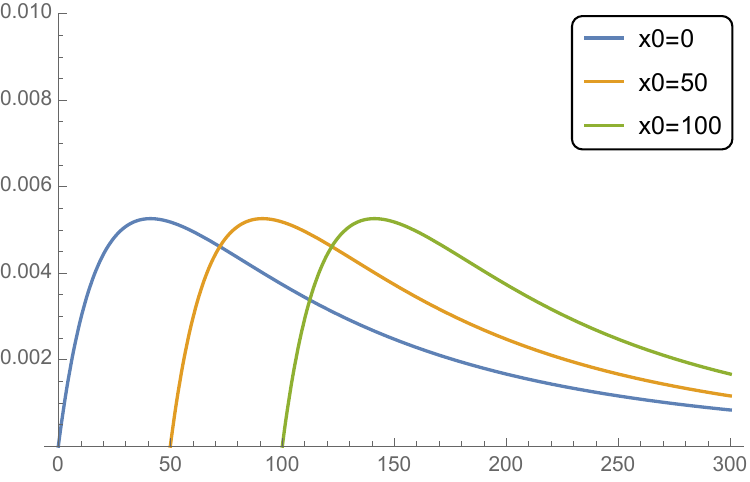}
	\end{subfigure}
	\caption{ \label{Parameters} Density plots of the IF distribution where, except for the parameter we vary in each plot, the parameters remain fixed to $p=1$, $b=1$, $c=200$, $q=2$ and $x_0=0$.}
\end{figure}

For the sake of illustration, we provide density plots of the IF distribution in Figure \ref{Parameters}. Except for the parameter we are varying, all the parameters remain fixed to $p=1$, $b=1$, $c=200$, $q=2$ and $x_0=0$.

Figure~\ref{Parameters} provides a visual inspection of the roles endorsed by the five parameters: $x_0\geqslant0$ is a location parameter (smaller than or equal to the lowest value of the data), $c>0$ a scale parameter, $q>0$ a tail-weight parameter, and $b\in\R_0$ a shape parameter regulating the skewness. By changing the sign of $b$ in IF, one gets  the  Inverse-IF distribution (such as, for instance, the Rayleigh and Inverse Rayleigh distribution, see below). A crucial role is played by $p\in[0,\infty]$ as it enables us to interpolate between power laws and power laws with exponential cut-off. We therefore name it \textsl{interpolation parameter}.


\section{Special cases and three main IF subfamilies} \label{sec:special}

One major appeal of the IF distribution is that it contains a plethora of well-known size distributions as special cases. For a clearer structure, we define three four-parameter subfamilies: 
\begin{itemize}
\item{the IF1 distribution where $p=0$,}
\item{the IF2 distribution where $p \to \infty$,}
\item{the IF3 distribution where $p\in(0,\infty)$ and $b=1$.}
\end{itemize}

Of course, there remain several other parameter combinations in the Interpolating Family, and perhaps in the future other interesting subfamilies will be given special attention. We  also by no means claim to be exhaustive in the special cases cited here, and the reader may well find other known distributions that are special cases of the IF but not mentioned here.


\subsection*{The IF1 distribution}

In the power law limit $p=0$, the pdf of the resulting four-parameter family of distributions, called \emph{Interpolating Family of the first kind~(IF1)}, is given by
$$
	f_0(x)=\text{sign}(b)q~g_0(x)G_0(x)^{-q-1}		=\frac{|b|q}{c} \left(\frac{x-x_0}{c}\right)^{b-1} \left(1+\left(\frac{x-x_0}{c}\right)^{b}\right)^{-q-1},
$$
where $x\in [x_0,\infty)$. Special cases of the IF1 distribution are, in decreasing order of the number of parameters,
the Lindsay--Burr type III distribution ($b<0$), the Pareto type IV distribution ($b>0$), the Dagum distribution ($b<0$ and $x_0=0$), the Pareto type II distribution ($b=1$), the Pareto type III distribution ($b>0$ and $q=1$), the Tadikamalla--Burr type XII distribution ($b>0$ and $x_0=0$), the Pareto type I distribution ($b=1$ and $c= x_0>0$), the Lomax distribution ($b=1$ and $x_0=0$), the Burr type XII distribution ($b>0, c=1$ and $x_0=0$) and the Fisk distribution ($b>0, q=1$ and $x_0=0$).

\subsection*{The IF2 distribution}

In the power law with exponential cut-off limit $p \to \infty$, the pdf of the resulting four-parameter family of distributions, called \emph{Interpolating Family of the second kind (IF2)}, is given by
$$f_{\infty}(x)= \text{sign}(b)q~g_{\infty}(x)G_{\infty}(x)^{-q-1} e^{-G_{\infty}(x)^{-q}}
			= \frac{|b|q}{c} \left(\frac{x-x_0}{c}\right)^{-bq-1} e^{-\left(\frac{x-x_0}{c}\right)^{-bq}},
$$
where $x\in[x_0,\infty)$. Special cases of the IF2 distribution are, in decreasing order of the number of parameters, the Weibull distribution ($b=-1$; if also $x_0=0$, we find the two-parameter Weibull distribution), the Fr\'{e}chet distribution ($b=1$; if also $x_0=0$, we find the two-parameter Fr\'{e}chet distribution), the Gumbel type II distribution ($b=1$ and $x_0=0$), the Rayleigh distribution ($b=-1,q=2$ and $x_0=0$), the Inverse Rayleigh distribution ($b=1, q=2$ and $x_0=0$), the Exponential distribution ($b=-1,q=1$ and $x_0=0$), and the Inverse Exponential distribution ($b=1,q=1$ and $x_0=0$).


\subsection*{The IF3 distribution}

The \emph{Interpolating Family of the third kind (IF3)} is characterized by  $0<p<\infty$ and $b=1$, resulting in the pdf
$$f_{p,1}(x) = \frac{q}{c} \left((p+1)^{-\frac{1}{q}} +\frac{x-x_0}{c}\right)^{-q-1} \left(1-\frac{1}{p+1}\left((p+1)^{-\frac{1}{q}} + \frac{x-x_0}{c}\right)^{-q}\right)^p,$$
where $x\in [x_0,\infty)$. 
Special cases   of the IF3 distribution are  the Generalized Lomax distribution ($x_0=0$) and the Stoppa distribution ($x_0=c(p+1)^{-\frac{1}{q}}$).


\subsection*{Distribution tree} 

A visual summary of the structure inherent to the IF distribution with its various special cases is given in Figure~\ref{fig:tree} below. Since the inverse of each distribution is obtained by switching the sign of the parameter $b$, we only give the tree for positive values of $b$. 

\begin{figure}
\begin{center}
\begin{sideways}
\begin{tikzpicture}
 \tikzstyle{block}=[rectangle, draw=blue!40, thick, fill=blue!10, text width=6em, text centered, rounded corners, minimum height=2em];
 \tikzstyle{indic}=[text width=6em];
 \tikzstyle{indic2}=[text width=15em];  
  \path (1.75,0) node [block] (IF) {IF}
            (-4.75,-3) node [block] (IF1) {IF1 \quad (Pareto IV)}
	   (1.75,-3) node [block] (F) {IF3}
            (8.25,-3) node [block] (IF2) {IF2}
            (-8,-6) node [block] (Par3) {Pareto III}
            (-4.75,-6) node [block] (TB12) {T-Burr XII}
            (-1.5,-6) node [block] (Par2) {Pareto II}
            (1.75,-6) node [block] (GL) {Generalized Lomax}
            (5,-6) node [block] (Stoppa) {Stoppa}
            (8.25,-6) node [block] (Frechet) {Fr\'{e}chet}
            (-8,-9) node [block] (Fisk) {Fisk}
            (-4.75,-9) node [block] (Burr3) {Burr XII}
            (-1.5,-9) node [block] (Lomax) {Lomax}
            (2.5,-9) node [block] (Par1) {Pareto I}
            (6.62,-9) node [block] (Gumbel) {Gumbel II}
            (5,-12) node [block] (IExponential) {Inverse Exponential}
            (8.25,-12) node [block] (IRayleigh) {Inverse Rayleigh}
            (-11,0) node [indic] (5p) {5 parameters}
            (-11,-3) node [indic] (4p) {4 parameters}
            (-11,-6) node [indic] (3p) {3 parameters}
            (-11,-9) node [indic] (2p) {2 parameters}
            (-11,-12) node [indic] (1p) {1 parameter};                            
  \draw[dashed] (IF) -- node[text width=4em, text centered]{$p=0$}(IF1);
  \draw[dashed] (IF) -- node[text width=4em, text centered]{$b=1$}(F);
  \draw[dashed] (IF) -- node[text width=4em, text centered]{$p\to \infty$}(IF2);
  \draw[dashed] (IF1) -- node[text width=4em, text centered]{$q=1$}(Par3);
  \draw[dashed] (IF1) -- node[text width=4em, text centered]{$x_0=0$}(TB12);
  \draw[dashed] (IF1) -- node[text width=4em, text centered]{$b=1$}(Par2);
  \draw[dashed] (F) -- node[text width=4em, text centered]{$p=0$}(Par2);
  \draw[dashed] (F) -- node[text width=4em, text centered]{$x_0=0$}(GL);
  \draw[dashed] (F) -- node[near end, text width=8em, text centered]{$x_0=c(p+1)^{-\frac{1}{q}}$}(Stoppa);
  \draw[dashed] (F) -- node[text width=4em, text centered]{$p\to \infty$}(Frechet); 
  \draw[dashed] (IF2) -- node[text width=4em, text centered]{$b=1$}(Frechet);
  \draw[dashed] (Par3) -- node[text width=4em, text centered]{$x_0=0$}(Fisk);
  \draw[dashed] (TB12) -- node[text width=4em, text centered]{$q=1$}(Fisk);
  \draw[dashed] (TB12) -- node[text width=4em, text centered]{$c=1$}(Burr3);
  \draw[dashed] (TB12) -- node[text width=4em, text centered]{$b=1$}(Lomax); 
  \draw[dashed] (Par2) -- node[text width=4em, text centered]{$x_0=0$}(Lomax); 
  \draw[dashed] (Par2) -- node[near end, text width=4em, text centered]{$x_0=c$}(Par1);
  \draw[dashed] (GL) -- node[near end, text width=4em, text centered]{$p=0$}(Lomax); 
  \draw[dashed] (Stoppa) -- node[near end, text width=4em, text centered]{$p=0$}(Par1);
  \draw[dashed] (GL) -- node[near end, text width=4em, text centered]{ $p \to \infty$}(Gumbel); 
  \draw[dashed] (Stoppa) -- node[text width=4em, text centered]{$p \to \infty$}(Gumbel);
  \draw[dashed] (Frechet) -- node[text width=4em, text centered]{$x_0=0$}(Gumbel); 
  \draw[dashed] (Gumbel) -- node[text width=4em, text centered]{$q=1$}(IExponential);
  \draw[dashed] (Gumbel) -- node[text width=4em, text centered]{$q=2$}(IRayleigh);                   
\end{tikzpicture}
\end{sideways}
\end{center}
\caption{Distribution tree of the Interpolating Family of size distributions} \label{fig:tree}
\end{figure}
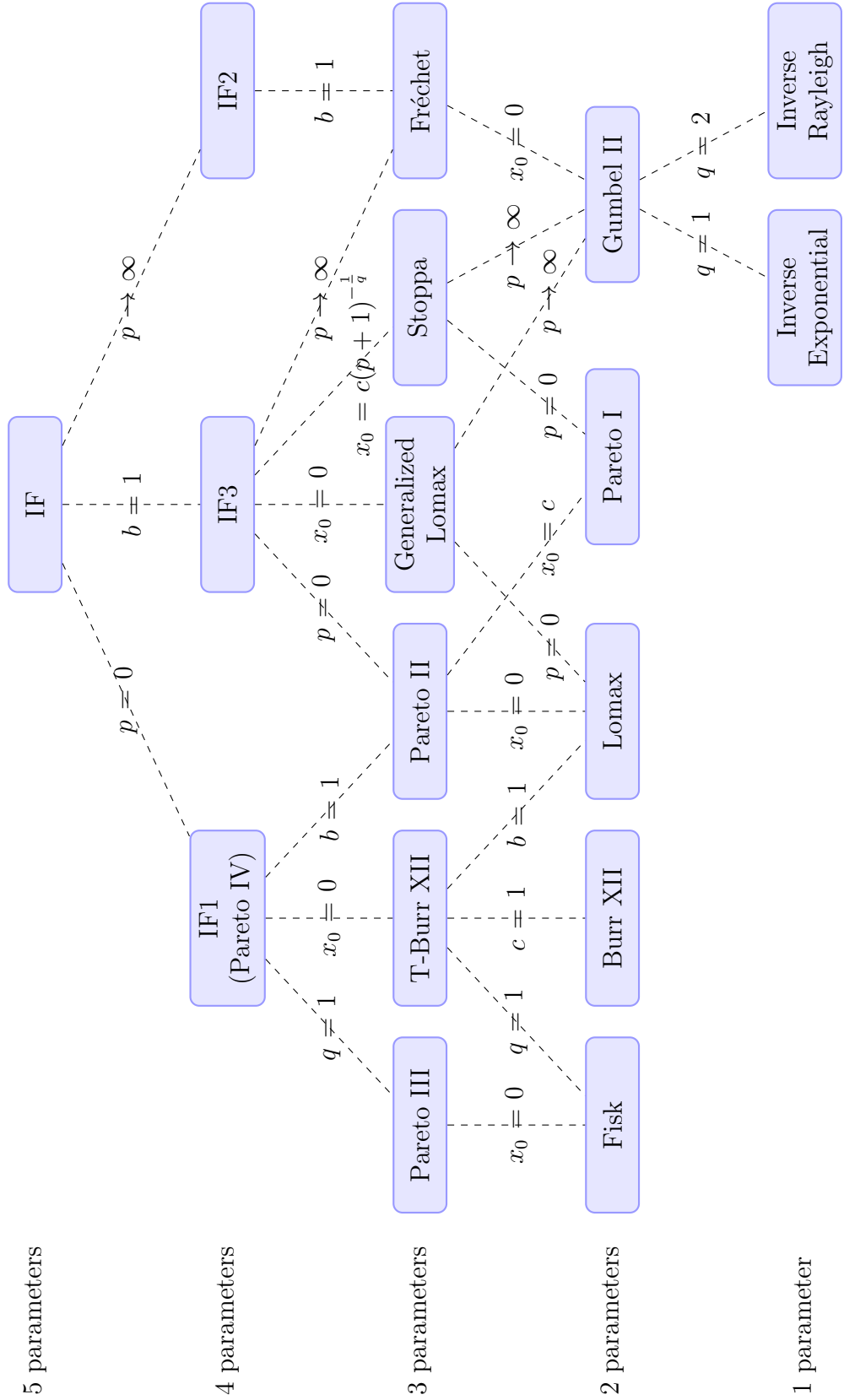


\section{Quantile-based properties} \label{sec:quant}

In this section we present and discuss quantile-based properties of the IF distribution. Since it contains so many special cases, the subsequent results provide in a single sweep those properties for the various distributions mentioned in Section~\ref{sec:special}.


\subsection{Cumulative distribution function, survival function and hazard function}
One major advantage of the IF distribution is that the cumulative distribution function (cdf) can be written under closed form:
\begin{align}
	F_p(x)
	 =\left \{ 
	 	\begin{array}{l l} 
			\left(1-\frac{1}{p+1}G_p(x)^{-q}\right)^{p+1} \quad &\text{if } b>0, \\
			1-\left(1-\frac{1}{p+1}G_p(x)^{-q}\right)^{p+1} \quad &\text{if } b<0.  \label{CDF}
		\end{array}
		\right.
\end{align}
Consequently, the survival or reliability function $S_p(x)=1-F_p(x)$ is extremely simple, too. The same holds true for the hazard function, defined as the ratio of the pdf and the survival function:
\begin{align*}
	H_p(x)=\frac{f_p(x)}{S_p(x)}
	= \left \{ 
	 	\begin{array}{l l} 
			q~g_p(x)G_p(x)^{-q-1}\frac{\left(1-\frac{1}{p+1}G_p(x)^{-q}\right)^{p}}{1-\left(1-\frac{1}{p+1}G_p(x)^{-q}\right)^{p+1}} \quad &\text{if } b>0, \\
			-q~g_p(x)G_p(x)^{-q-1} \frac{1}{\left(1-\frac{1}{p+1}G_p(x)^{-q}\right)}\quad &\text{if } b<0 .
		\end{array}
		\right.
\end{align*}


\subsection{Quantile function and median}

Very conveniently, the quantile function takes a nice form thanks to the simple expression of the cdf~\eqref{CDF}. Given the wide range of quantile-based statistical tools and methods such as QQ-plots, interquartile range or quantile regression, this is a very welcomed feature of the IF distribution.
For $b>0$, the quantile function is given by
\begin{align*}
	Q_p^+(y)=F_p^{-1}(y) = \left\{
				\begin{array}{l l}
					x_0+c\left((1-y)^{-\frac{1}{q}}-1\right)^{\frac{1}{b}} &\text{if } p=0,\\
					x_0+c(p+1)^{-\frac{1}{bq}}\left(\left(1-y^{\frac{1}{p+1}}\right)^{-\frac{1}{q}}-1\right)^{\frac{1}{b}} &\text{if }  0< p<\infty, \\
					x_0+c\left(\ln \left(\frac{1}{y}\right)\right)^{-\frac{1}{bq}} &\text{if } p \to \infty,
				\end{array}
			\right.
\end{align*}
for $y\in[0,1]$. The expression for $b<0$ is readily obtained via the relationship $Q_p^-(y)=Q_p^+(1-y)$, and we define the quantile function $Q_p(y)$ as $Q_p^+(y)$ if $b>0$ and as $Q_p^-(y)$ if $b<0$. The median is  uniquely defined as 
\begin{align*}
	\mbox{Median}= \left\{
				\begin{array}{l l}
					x_0+c\left((2^{\frac{1}{q}}-1\right)^{\frac{1}{b}} &\text{if } p=0,\\
					x_0+c(p+1)^{-\frac{1}{bq}}\left(\left(1-2^{-\frac{1}{p+1}}\right)^{-\frac{1}{q}}-1\right)^{\frac{1}{b}} &\text{if }  0< p<\infty, \\
					x_0+c \left( \ln \left(2\right)\right)^{-\frac{1}{bq}} &\text{if } p \to \infty.
				\end{array}
			\right.
\end{align*}


\subsection{Random variable generation}
The closed form of the quantile functions entails a straightforward  random variable generation process from the IF. Indeed, it suffices to generate a random variable $U$ from a uniform distribution on the interval $[0,1]$, and then apply $Q_p$ to it. The resulting random variable $Q_p(U)$ follows the IF~distribution. The simplicity of the procedure is particularly important for Monte Carlo simulation purposes and shows that all size distributions that are part of the Interpolating Family  benefit from a straightforward random variable generation procedure.


\section{Moments, mean and variance}\label{sec:moments}

We now provide the general moment expressions of the IF  distribution. Particular focus shall be given to the mean and variance expressions, which we can write out explicitly in terms of the gamma and beta functions. We conclude the section with a table containing the mean expressions for the various size distributions mentioned in Section~\ref{sec:special}.

\hfill

The $r^\text{th}$ moment of the IF distribution is given by 
\begin{equation*}
	\begin{array}{ccl}
		\mathbb{E}\left[X^r\right]&=& \int \limits_{x_0}^\infty x^r~\text{IF}(x;p,b,c,q,x_0)~\mathrm{d}x.
	\end{array}
\end{equation*}
We will first treat the cases when $p<\infty$. Making the change of variables $y= (p+1)^{-\frac{1}{q}} + \left(\frac{x-x_0}{c}\right)^b$ and applying Newton's binomial theorem,  we get
\begin{equation*}
	\begin{array}{ccl}
		\mathbb{E}\left[X^r\right]&=& \sum \limits_{i=0}^r {r \choose i} x_0^i c^{r-i} \underbrace{ \int \limits_{(p+1)^{-\frac{1}{q}}}^\infty q~y^{-q-1} \left( y-(p+1)^{-\frac{1}{q}} \right)^{\frac{r-i}{b}} \left( 1 - \frac{y^{-q}}{p+1} \right)^p \mathrm{d}y.}_{I(p,b,q)}
	\end{array}
\end{equation*}
Note how the sign of $b$ vanishes during this change of variable. The following result establishes under which conditions the moments of the IF distribution exist and are finite. 

\begin{prop}\label{propmom}
The $r^\text{th}$ moment of the IF distribution for $p<\infty$ exists and is finite  if and only if  $b>0$ and $r<bq$ or $b<0$ and $r<-b$. 
\end{prop}

The proof is provided in the Appendix.  It is in principle possible to write out $I(p,b,q)$ as an infinite series of beta functions, but since this expression is rather intricate and  needs to be worked out on a case-by-case basis just like $I(p,b,q)$, we refrain from doing so. However, in what follows, we compute the integral $I(p,b,q)$ for the four-parameter distributions IF1 and~IF3 and obtain there more tractable expressions.

\hfill

\subsection*{Moments of the IF1 distribution}
If we plug in $p=0$ and then set $z=\frac{1}{y}$ we get
\begin{equation*}
	\begin{array}{ccl}
		 I(0,b,q) &=& \int \limits_1^\infty q~y^{-q-1} \left( y-1 \right)^{\frac{r-i}{b}} \mathrm{d}y = q \int \limits_0^1 z^{q-1-\frac{r-i}{b}} \left(1-z\right)^{\frac{r-i}{b}} \mathrm{d}z.
	\end{array}
\end{equation*}	
If either $b>0$ and $r<bq$ or else $b<0$ and $r<-b$ (i.e., under the conditions identified in Proposition~\ref{propmom}), then this can be written as
\begin{equation*}
	\begin{array}{ccl}
		I(p,b,q) &=& q B(q-\frac{r-i}{b},1+ \frac{r-i}{b}), 
	\end{array}
\end{equation*}
where $B(\cdot,\cdot)$ stands for the beta function. The $r^\text{th}$ moment of the IF1 distribution is thus given by
\begin{equation*}
	\begin{array}{ccl}
		\mathbb{E}\left[X^r\right]&=& \sum \limits_{i=0}^r {r \choose i} x_0^i c^{r-i} q B(q-\frac{r-i}{b},1+ \frac{r-i}{b}) \quad \text{if } \left\{
\begin{array}{l}
b>0 \text{ and } r<bq,\\
b<0 \text{ and } r<-b,
\end{array}
\right.
	\end{array}\vspace{2mm}
\end{equation*}
otherwise the moment does not exist. The mean and variance of the IF1 distribution  are then respectively given by
	\begin{equation}\label{mean}\mathbb{E}(X)=x_0+cq~B\left(q-\frac{1}{b},1+\frac{1}{b}\right) \quad \text{if } \left \{ \begin{array}{ll}
																	b>\frac{1}{q}, \\
																	b<-1,\\
																	\end{array}
																	\right.\end{equation}
	and																
	$$\mathbb{V}(X)=c^2\left[q~B\left(q-\frac{2}{b},1+\frac{2}{b}\right)-\left(q~B\left(q-\frac{1}{b},1+\frac{1}{b}\right)\right)^2\right]\quad \text{if } \left \{ \begin{array}{ll}
																	b>\frac{2}{q}, \\
																	b<-2.\\
																	\end{array}\right.$$
Note that sometimes it can be convenient to rewrite the mean~\eqref{mean} under the form 
$$
\mathbb{E}(X)=x_0+c\frac{\Gamma\left(q-\frac{1}{b}\right)\Gamma\left(1+\frac{1}{b}\right)}{\Gamma(q)}
$$
with $\Gamma(\cdot)$ the gamma function. Special cases of the mean expressions can be found in Table~\ref{tab1}.

\subsection*{Moments of the IF3 distribution}
If $b=1$, for $r<q$ we make the change of variables $z=\frac{y^{-q}}{p+1}$ and get
\begin{equation*}
	\begin{array}{ccl}
		I(p,1,q) &=& (p+1)^{1-\frac{r-i}{q}} \int \limits_0^1 \left(z^{-\frac{1}{q}}-1\right)^{r-i} \left(1-z\right)^p \mathrm{d}z \\
		&=& (p+1)^{1-\frac{r-i}{q}} \sum \limits_{k=0}^{r-i} {r-i \choose k} (-1)^k \int \limits_0^1 z^{-\frac{1}{q}(r-i-k)}\left(1-z\right)^p\mathrm{d}z\\
		&=& (p+1)^{1-\frac{r-i}{q}} \sum \limits_{k=0}^{r-i} {r-i \choose k} (-1)^k B\left(1-\frac{1}{q}(r-i-k),p+1\right).
	\end{array}
\end{equation*}
These manipulations are possible since the integral is finite under $r<q$. Consequently, the $r^\text{th}$ moment of the IF3 distribution  is given by
\begin{equation*}
	\begin{array}{ccl}
		\mathbb{E}\left[X^r\right]&=& \sum \limits_{i=0}^r {r \choose i} x_0^i c^{r-i}~(p+1)^{1-\frac{r-i}{q}} \sum \limits_{k=0}^{r-i} {r-i \choose k} (-1)^k B\left(1-\frac{1}{q}(r-i-k),p+1\right).
	\end{array}
\end{equation*}
The associated mean and variance  are
	$$\mathbb{E}(X)=x_0+c (p+1)^{1-\frac{1}{q}} \left(B\left(1-\frac{1}{q},p+1\right)-\frac{1}{p+1}\right) \quad \text{if } 1<q,$$
	and
	\begin{eqnarray*}
	\mathbb{V}(X)&=&c^2(p+1)^{1-\frac{2}{q}}\left[B\left(1-\frac{2}{q},p+1\right)-2B\left(1-\frac{1}{q},p+1\right)+\frac{1}{p+1}\right]\\
	&&-c^2(p+1)^{2-\frac{2}{q}}\left[B\left(1-\frac{1}{q},p+1\right)-\frac{1}{p+1}\right]^2 \quad \text{if } 2<q.		\end{eqnarray*}
The special cases of the mean expressions for the Generalized Lomax and Stoppa distributions can be found in Table~\ref{tab1}.

\

Let us now consider the case $p=\infty$.

\subsection*{Moments of the IF2 distribution}

The $r$-th moment is calculated as
\begin{equation*}
		\mathbb{E}\left[X^r\right]= \int \limits_{x_0}^\infty x^r~\frac{|b|q}{c} \left(\frac{x-x_0}{c}\right)^{-bq-1} e^{-\left(\frac{x-x_0}{c}\right)^{-bq}}~\mathrm{d}x.
\end{equation*}
In this case the finite moments conditions are more easily seen and require no formal statement under the form of a proposition. When $b<0$, all moments exist, while for $b>0$ we can see that the integrand behaves like $x^{r-bq-1}$ for large values of $x$, implying existence of thr $r$-th moment iff $r<bq$. Under these conditions, the change of variables $y=  \left(\frac{x-x_0}{c}\right)^{-bq}$ combined with Newton's binomial theorem implies that the integral  is equal to
\begin{equation*}
		\sum_{i=0}^r\binom{r}{i}x_0^ic^{r-i}\int \limits_{0}^\infty y^{-\frac{r-i}{bq}} e^{-y}~\mathrm{d}y=\sum_{i=0}^r\binom{r}{i}x_0^ic^{r-i}\Gamma\left(1-\frac{r-i}{bq}\right)
\end{equation*}
by definition of the gamma function.
The $r^\text{th}$ moment of the IF2 distribution is therefore given by
\begin{equation*}
	\begin{array}{ccl}
		\mathbb{E}\left[X^r\right]&=& \sum \limits_{i=0}^r {r \choose i} x_0^i c^{r-i}~\Gamma\left(1- \frac{r-i}{bq}\right) \quad \text{if } \left\{
		\begin{array}{l}
b>0 \text{ and } r<bq,\\
b<0,
\end{array}
\right.
	\end{array}
\end{equation*}
otherwise the moment does not exist. The corresponding mean and variance take on the expressions
	$$\mathbb{E}(X)=x_0+c~\Gamma\left(1-\frac{1}{bq}\right) \quad \text{if } \left \{ \begin{array}{ll}
																	b>\frac{1}{q}, \\
																	b<0,\\
																	\end{array}
																	\right.$$
	and																
	$$\mathbb{V}(X)=c^2 \left[\Gamma\left(1-\frac{2}{bq}\right)-\left(\Gamma\left(1-\frac{1}{bq}\right)\right)^2\right] \quad \text{if } \left \{ \begin{array}{ll}
																	b> \frac{2}{q}, \\
																	b<0.&\\
																	\end{array}
																	\right.$$																Special cases of the mean expressions can be found in Table~\ref{tab1}.

\begin{table}\small
\begin{center}
	\begin{tabular}{|c||c||c|c|c|c|}
	\hline
	Distribution & \# & Parameters & Mean & Constraint \\
	name & par. & $(p,b,c,q,x_0)$ & $\mathbb{E}\left[X\right]$ & \\
	\hline
	\hline
	Pareto IV & 4 &$(0,\frac{1}{\gamma}>0,c,q,x_0)$ & $x_0+cqB\left(q-\gamma,1+\gamma \right)$& $\gamma < q$ \\
	\hline
	Lindsay--Burr III & 4 &$(0,b<0,c,q,x_0)$ & $x_0+cqB\left(q-\frac{1}{b},1+\frac{1}{b}\right)$& $b<-1$ \\
	\hline
	Dagum  & 3 &$(0,b<0,c,q,0)$ & $cqB\left(q-\frac{1}{b},1+\frac{1}{b}\right)$& $b<-1$ \\
	\hline
	Pareto II & 3 & $(0,1,c,q,x_0)$ & $x_0+\frac{c}{q-1}$ & $q>1$\\
	\hline
	Pareto III & 3 & $(0,\frac{1}{\gamma}>0,c,1,x_0)$  & $x_0 + c \Gamma\left(1-\gamma \right)\Gamma\left(1+\gamma \right)$& $\gamma < 1$\\
	\hline 
	Tadikamalla--Burr XII & 3 & $(0,b>0,c,q,0)$ & $cq B\left(q-\frac{1}{b},1+\frac{1}{b}\right)$ & $bq>1$  \\
	\hline
	Fisk & 2 & $(0,b>0,c,1,0)$ &  $c \Gamma\left(1-\frac{1}{b}\right) \Gamma\left(1+\frac{1}{b}\right)$ &$b>1$ \\
	\hline
	Lomax & 2 & $(0,1,c,q,0)$ & $\frac{c}{q-1}$ & $q>1$\\
	\hline
	Pareto I & 2 & $(0,1,x_0,q,x_0)$ & $\frac{q}{q-1}x_0$& $q>1$\\
	\hline
	Burr XII  & 2 & $(0,b>0,1,q,0)$& $qB\left(q-\frac{1}{b},1+\frac{1}{b}\right) $& $bq>1$  \\
	\hline
	\hline
	Weibull & 3 & $(\infty,-1,c,q,x_0)$&$x_0+c\Gamma\left(1+\frac{1}{q}\right)$ & \\
	\hline
	Fr\'{e}chet & 3 & $(\infty,1,c,q,x_0)$ & $x_0+ c\Gamma\left(1-\frac{1}{q}\right)$ & $q>1$\\
	\hline
	Gumbel II & 2 & $(\infty,1,c,q,0)$ & $c\Gamma\left(1-\frac{1}{q}\right)$ & $q>1$\\
	\hline
	Rayleigh & 1 &$(\infty,-1,c,2,0)$ & $\frac{c}{2} \sqrt{\pi}$ & \\
	\hline
	Inverse Rayleigh & 1 &$(\infty,1,c,2,0)$ & $c \sqrt{\pi}$ & \\
	\hline
Exponential & 1 & $(\infty,-1,c,1,0)$& $c$ & \\
	\hline
	Inverse Exponential & 1 & $(\infty,1,c,1,0)$& Not defined & Violated\\
	\hline
	\hline
	Generalized Lomax & 3 &$(m-1,1,c,q,0)$& $c~m^{1-\frac{1}{q}} \left(B\left(1-\frac{1}{q},m\right)-\frac{1}{m}\right)$ &$q>1$ \\
	\hline
	Stoppa & 3 & $(m-1,1,c,q,cm^{-\frac{1}{q}})$ & $c~m^{1-\frac{1}{q}}~B\left(1-\frac{1}{q},m\right)$&$q>1$ \\
	\hline
	\end{tabular}
\end{center}
\caption{Expressions for the mean} \label{tab1}
\end{table}
\normalsize

 \section{Unimodality and location of the mode}\label{sec:mode}

Determining the mode of a distribution is an important issue, which we tackle in this section.  We study the derivative of $x\mapsto f_p(x)$, with particular emphasis on the three main subfamilies IF1, IF2 and~IF3 described in Section~\ref{sec:special}.
As we show in the Appendix, the derivative of the pdf vanishes either at the boundary $x=x_0$ of the domain or at 
\begin{equation*}
x= x_0+c (p+1)^{-\frac{1}{bq}} \left(t^{-\frac{1}{q}}-1\right)^{\frac{1}{b}},
\end{equation*}
where $t$ is solution of the almost cyclic equation
\begin{equation}\label{cyclicequation}
(b-1) t^{-\frac{1}{q}} (1-t) -b(q+1)(t^{-\frac{1}{q}}-1)(1-t) + pbq (t^{-\frac{1}{q}}-1) t =0.
\end{equation}
This allows us to draw the following conclusions regarding the modality of the IF distribution.

\begin{itemize} 
\item The mode of the IF1 distribution ($p=0$) is given by
\begin{equation*}
\left\{
\begin{array}{ll}
x_0 & \text{ if }b = -\frac{1}{q} \text{ or } b=1,\\
x_0+c\left(\frac{b-1}{bq+1}\right)^{\frac{1}{b}} & \text{ if } b<-\frac{1}{q} \text{ or } b>1,
\end{array}
\right.
\end{equation*}
whereas in the remaining cases, i.e.~$b \in ]-\frac{1}{q};1[$, there is a vertical asymptote at $x=x_0$. We plot in Figure~\ref{fig1} a contour plot of the mode of the IF1 distribution.
\begin{figure}[h!]
	\begin{center}
	\includegraphics[scale=0.7]{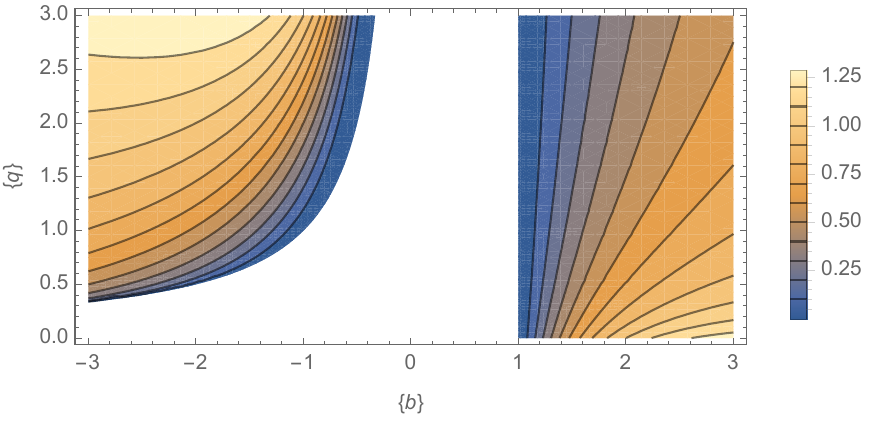}
	\caption{Contour plot of the mode of the IF1 distribution as a function of the parameters, with $c=1$ and $x_0=0$.}\label{fig1}
	\end{center}  
\end{figure}
\item The mode of the IF2 distribution ($p\to \infty$) is given by
\begin{equation*}
\left\{
\begin{array}{ll}
x_0 & \text{ if }b = -\frac{1}{q},\\
x_0+c\left(\frac{bq}{bq+1}\right)^{\frac{1}{bq}} & \text{ if } b<-\frac{1}{q} \text{ or } b>0,
\end{array}
\right.
\end{equation*}
whereas in the remaining cases, i.e.~$b \in ]-\frac{1}{q},0]$, there is a vertical asymptote at $x=x_0$. Figure~\ref{fig2} shows a contour plot of the mode of the IF2 distribution.
\begin{figure}[h!]
	\begin{center}
	\includegraphics[scale=0.7]{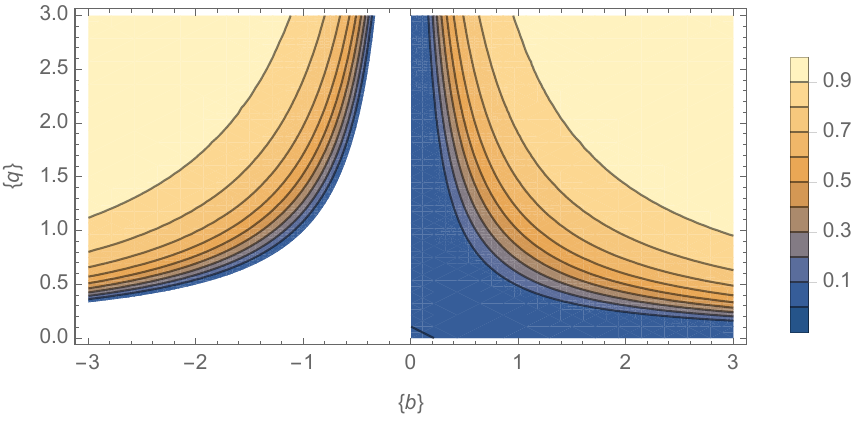}
	\caption{Contour plot of the mode of the IF2 distribution as a function of the parameters, with $c=1$ and $x_0=0$.}\label{fig2}
	\end{center}  
\end{figure}
\item The mode of the IF3 distribution ($0<p<\infty$ and $b=1$) is given by
$$
x_0+c(p+1)^{-\frac{1}{q}}\left(\left(\frac{q+1}{(p+1)q+1}\right)^{-\frac{1}{q}}-1\right).
$$
A contour plot of the mode of the IF3 distribution can be seen in Figure~\ref{fig3}.
\begin{figure}[h!]
	\begin{center}
	\includegraphics[scale=0.7]{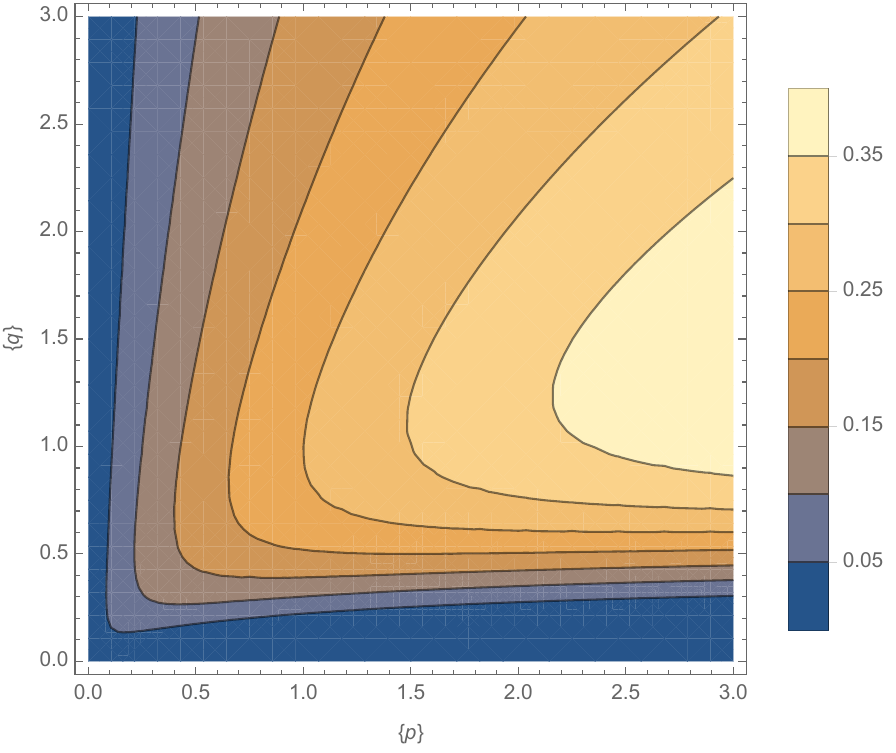}
	\caption{Contour plot of the mode of the IF3 distribution as a function of the parameters, with $c=1$ and $x_0=0$.}\label{fig3}
	\end{center}  
\end{figure}
\end{itemize}

For calculation details, see the Appendix. We note that all three subfamilies are unimodal, which is coherent with the related special cases from the literature. Moreover, we have derived the exact expressions of the modes. This unimodality is a very attractive feature from an interpretability point of view: bimodal or multimodal distributions are arguably best modeled as a mixture of unimodal distributions. It is therefore not surprising that many new distributions are built with the target of being unimodal; for modern examples, see e.g.~\cite{jones2014generating, katojones2015, fujisawaabe2015}.

\section{Final remarks}\label{sec:final}

We have built in this paper an overarching family of size distributions, the Interpolating Family of distributions, and shown how it indeed interpolates between power laws and power laws with exponential cut-off. This sheds interesting new light on these highly distinct types of size distributions, and we hope that our construction inspired from statistical physics will stimulate researchers to search for bridges between other apparently unrelated distributions. Understanding the links between probability laws and grouping them into classes with similar properties has become particularly important nowadays, given the plethora of new distributions. We refer the interested reader to the review papers \cite{Jones2015} and \cite{babic2019comparison} for further information about classifying flexible distributions for univariate respectively multivariate data, and to \cite{ley2020flexible}  for an overview and discussion of advantages and limitations of flexible models.


Finally, we recall that the aim of this research was to treat power laws and power laws with exponential cut-off in a unified way, and to develop general properties for the  IF distribution. We purposely did not provide inferential procedures for the full 5-parameter IF distribution since we have noticed that distinct combinations of the five parameters lead to nearly identical shapes of the density, and consequently maximum likelihood estimation may be ill-defined. We therefore recommend to restrict to subfamilies for inferential purposes\footnote{Given the simple form of the cumulative distribution function, special cases of the IF  are tailor-made for dealing with censored data.}, or to consider interesting new special cases of  the IF distribution. Our results readily yield the theoretical basis for such future research.

\vspace{1cm}

\section*{Appendix}

\subsection*{Proof of Proposition~\ref{propmom}}
 We start by performing the change of variables $y^{-q}/(p+1)=t$ inside the integral $I(p,b,q)$, yielding
$$
(p+1)^{1-\frac{r-i}{qb}}\int_0^1\left(t^{-1/q}-1\right)^{\frac{r-i}{b}}(1-t)^p\,dt.
$$
Since $p<\infty$ and $t\in[0,1]$, the factor $(1-t)^p$ is always bounded by 1 and hence causes no problem. Concentrating our attention on $\left(t^{-1/q}-1\right)^{\frac{r-i}{b}}$, we need to distinguish two cases:
\begin{itemize}
\item[$b>0$:] The finiteness of the integral depends on  $\left(t^{-1/q}-1\right)^{\frac{r-i}{b}}$ when $t$ approaches 0, in which case the expression behaves like $t^{\frac{i-r}{qb}}$ and hence is finite iff $\frac{i-r}{qb}>-1$, which is equivalent to $r-i<qb$. Since this needs to hold for every $i\in[0,r]$, we conclude that the $r$-th moment is finite iff $r<bq$.
\item[$b<0$:]  The finiteness of the integral depends on  $\left(t^{-1/q}-1\right)^{\frac{r-i}{b}}$ when $t$ approaches 1, in which case the expression behaves like $(1-t)^{\frac{r-i}{b}}$ and hence is finite iff $\frac{r-i}{b}>-1$, which is equivalent to $r-i<-b$. Since this needs to hold for every $i\in[0,r]$, we conclude that the $r$-th moment is finite iff $r<-b$.
\end{itemize}
This concludes the proof. \hfill $\square$

\subsection*{Mode calculation}

The derivative of the pdf~(\ref{IF}) vanishes if and only if \small
\begin{eqnarray*}
	0=&(b-1) \left(\frac{x-x_0}{c}\right)^{b-2} \left((p+1)^{-\frac{1}{q}}+\left(\frac{x-x_0}{c}\right)^b\right)^{-q-1} \left(1-\frac{1}{p+1}\left( (p+1)^{-\frac{1}{q}}+\left(\frac{x-x_0}{c}\right)^b\right)^{-q}\right)^p \\
	&-b(q+1)\left(\frac{x-x_0}{c}\right)^{2b-2} \left((p+1)^{-\frac{1}{q}}+\left(\frac{x-x_0}{c}\right)^b\right)^{-q-2} \left(1-\frac{1}{p+1}\left( (p+1)^{-\frac{1}{q}}+\left(\frac{x-x_0}{c}\right)^b\right)^{-q}\right)^p \\
	&+\frac{pbq}{p+1} \left(\frac{x-x_0}{c}\right)^{2b-2} \left((p+1)^{-\frac{1}{q}}+\left(\frac{x-x_0}{c}\right)^b\right)^{-2q-2}\left(1-\frac{1}{p+1}\left( (p+1)^{-\frac{1}{q}}+\left(\frac{x-x_0}{c}\right)^b\right)^{-q}\right)^{p-1}. 
\end{eqnarray*} \normalsize
If we set $y= \frac{x-x_0}{c}$, then the above holds true if either 
\small \begin{equation*}
y^{b-2}
\left((p+1)^{-\frac{1}{q}}+y^b\right)^{-q-2} 
\left(1-\frac{1}{p+1}\left( (p+1)^{-\frac{1}{q}}+y^b\right)^{-q}\right)^{p-1}=0
\end{equation*} \normalsize
or
\small \begin{eqnarray}\label{secondequation}
&&0 = (b-1) \left((p+1)^{-\frac{1}{q}}+y^b\right) \left(1-\frac{1}{p+1}\left( (p+1)^{-\frac{1}{q}}+y^b\right)^{-q}\right) \\
& & \quad\quad -b (q+1) y^{b} \left(1-\frac{1}{p+1}\left( (p+1)^{-\frac{1}{q}}+y^b\right)^{-q}\right)  + \frac{pbq}{p+1} y^{b} \left((p+1)^{-\frac{1}{q}}+y^b\right)^{-q}. \nonumber
\end{eqnarray} \normalsize
The only solution that the first equation can possibly admit is $y=0$. This corresponds to $x=x_0$, i.e.~to the boundary of the domain. On the other hand, equation~(\ref{secondequation}) may admit interior solutions. We separate the analysis of equation~(\ref{secondequation}) in two parts: first the case $p$ finite from which we deduce the mode of the IF1 and IF3 subfamilies and second the case $p \to \infty$ which gives the mode of the IF2~subfamily. If $p$ is finite, then we set 
$$t=\frac{1}{p+1}\left((p+1)^{-\frac{1}{q}}+y^b\right)^{-q}$$ and equation~(\ref{secondequation}) simplifies to the almost cyclic equation~(\ref{cyclicequation}):
$$(b-1) t^{-\frac{1}{q}}(1-t)-b(q+1)\left(t^{-\frac{1}{q}}-1\right)(1-t)+pbq\left(t^{-\frac{1}{q}}-1\right)t=0.$$
Solving this equation in all generality is possible numerically but we will restrict ourselves to show how to get closed-form solutions for the two subfamilies IF1 and~IF3.
For the IF1 distribution ($p=0$), equation~(\ref{cyclicequation}) further simplifies to
\begin{align*}
	&(b-1)t^{-\frac{1}{q}} (1-t)-b(q+1)\left(t^{-\frac{1}{q}}-1\right)(1-t) =0. 
\end{align*}
While we recover the boundary solution $x=x_0$ if $b>0$, we also find an interior solution $x=x_0+c \left(\frac{b-1}{bq+1}\right)^{\frac{1}{b}}$ if either $b<-\frac{1}{q}$ or $b>1$. Repeating the procedure with the second derivative of the pdf~(\ref{IF}), a straightforward but tedious calculation shows that the interior solution thus found indeed corresponds to a maximum and that the mode occurs on the boundary $x=x_0$ if either $b=-\frac{1}{q}$ or $b=1$.

For the IF3 distribution ($0<p<\infty$ and $b=1$), equation~(\ref{cyclicequation}) further simplifies to
\begin{align*}
	&-(q+1)\left(t^{-\frac{1}{q}}-1\right)(1-t)+pq\left(t^{-\frac{1}{q}}-1\right)t =0. 
\end{align*}
This equation admits two solutions: the boundary solution $x=x_0$ and the interior solution $x=x_0+c (p+1)^{-\frac{1}{q}} \left(\left(\frac{q+1}{(p+1)q+1}\right)^{-\frac{1}{q}}-1\right)$. One can then check that the latter corresponds to the mode of the IF3~distribution and that this mode, and thus the interior solution, moves towards the boundary as $p$ and $q$ tend to zero.

On the other hand, for the IF2 distribution ($p\to \infty$), equation~(\ref{secondequation}) simplifies to
$$0 = (b-1) y^b - b(q+1)y^b + bq  y^{b-bq}.$$
We deduce that the derivative of the pdf of the IF2 vanishes either at the boundary $x=x_0$ if $b>0$ or at the interior point $x=x_0+c\left(\frac{bq}{bq+1}\right)^{\frac{1}{bq}}$ if $b<-\frac{1}{q}$ or $b>0$. Similarly as for the IF1, tedious second derivative calculations reveal that the interior solution always corresponds to a maximum and that the mode occurs on the boundary if either $b=-\frac{1}{q}$ or $b=0$.

\bibliographystyle{apalike}
\bibliography{IFrefs}

\end{document}